\documentclass[final]{siamltex}
\usepackage{epsfig}
\usepackage{amssymb}
\usepackage{amsmath}
\newtheorem{conj}[theorem]{Conjecture}

\newcounter{Examplecount}
\newenvironment{example}[1][Example \arabic{Examplecount}.]{\begin{trivlist}
\item[\hskip \labelsep {\bfseries #1}]}{\end{trivlist}\stepcounter{Examplecount}}

\newcommand\beq{\begin{equation}}
\newcommand\eeq{\end{equation}}
\newcommand\bce{\begin{center}}
\newcommand\ece{\end{center}}
\newcommand\bea{\begin{eqnarray}}
\newcommand\eea{\end{eqnarray}}
\newcommand\ben{\begin{enumerate}}
\newcommand\een{\end{enumerate}}
\newcommand\bit{\begin{itemize}}
\newcommand\eit{\end{itemize}}
\newcommand\brr{\begin{array}}
\newcommand\err{\end{array}}
\newcommand\bt{\begin{tabular}}
\newcommand\et{\end{tabular}}

\newcommand\nn{\nonumber}
\newcommand\bs{\bigskip}
\newcommand\ms{\medskip}
\newcommand\ul{\underline}

\renewcommand\S{{\mathcal S}}
\renewcommand\H{{\mathcal H}}
\newcommand\T{{\mathcal T}}
\newcommand\E{{\mathcal E}}
\newcommand\W{{\mathcal W}}
\newcommand\Al{\operatorname{Allow}}
\newcommand\F{\operatorname{Forb}}
\newcommand\B{\operatorname{B}}
\newcommand\A{\operatorname{Av}}

\newcommand\Pat{\operatorname{Pat}}
\newcommand\des{\operatorname{des}}
\newcommand\hp{\hat\pi}

\author{Sergi Elizalde \thanks{Department of Mathematics, Dartmouth College, Hanover, NH 03755-3551 ({\tt sergi.elizalde@dartmouth.edu})}}
\title{The number of permutations realized by a shift}

\begin{document}
\maketitle

\begin{abstract}
A permutation $\pi$ is realized by the shift on $N$ symbols if there is
an infinite word on an $N$-letter alphabet whose successive left shifts by one position are lexicographically in the same relative order as $\pi$.
The set of realized permutations is closed under consecutive pattern containment. Permutations that cannot be realized are called forbidden patterns.
It was shown in~\cite{AEK} that the shortest forbidden patterns of the shift on $N$ symbols have length $N+2$.
In this paper we give a characterization of the set of permutations that are realized by the shift on $N$ symbols, and we enumerate them according to their length.
\end{abstract}

\begin{keywords}
shift, consecutive pattern, forbidden pattern
\end{keywords}

\begin{AMS}
Primary 05A05; Secondary 05A15, 37M10, 94A55
\end{AMS}

\section{Introduction and definitions}

The original motivation for this paper comes from an innovative application of pattern-avoiding permutations to dynamical systems (see~\cite{AEK,Bandt}), which is based on the following idea.
Given a piecewise monotone map on a one-dimensional interval, consider the finite sequences (orbits) that are obtained by iterating the map,
starting from any point in the interval. It turns out that the relative order of the entries in these sequences cannot be arbitrary.
This means that, for any given such map, there will be some order patterns that will never appear in any orbit. The set of such patterns, which we call forbidden patterns,
is closed under consecutive pattern containment. These facts can be used to distinguish random from deterministic time series.

A natural question that arises is how to determine, for a given map, what its forbidden patterns are. While this problem is wide open in general, in the present paper
we study it for a particular kind of maps, called (one-sided) shift systems. Shift systems are interesting for two reasons.
On one hand, they exhibit all important features of low-dimensional chaos, such as sensitivity to initial conditions, strong mixing, and a dense set of periodic
points. On the other hand, they are natural maps from a combinatorial perspective, and the study of their forbidden patterns can be done in an elegant combinatorial way.

Forbidden patterns in shift systems were first considered in~\cite{AEK}. The authors prove that the smallest forbidden pattern of the shift on $N$ symbols has length $N+2$. They also conjecture
that, for any $N$, there are exactly six forbidden patterns of minimal length.
In the present paper we give a complete characterization of forbidden patterns of shift systems, and enumerate them with respect to their length.

We will start by defining consecutive pattern containment, forbidden patterns in maps, and shift systems, and introducing some notations and background.
In Section~\ref{sec:numsymb} we give a formula for the parameter that determines how many symbols are needed in order for a permutation to be realized by a shift.
This characterizes allowed and forbidden patterns of shift maps. In Section~\ref{sec:mostsymb} we give another equivalent characterization involving a transformation on permutations, and we prove
that the shift on $N$ symbols has six forbidden patterns of minimal length $N+2$, as conjectured in~\cite{AEK}.
In Section~\ref{sec:binary} we give a formula for the number of patterns of a given length that are realized by the binary shift. In Section~\ref{sec:anyN} we generalize the results from
Section~\ref{sec:binary} to the shift on $N$ symbols, for arbitrary $N$. We end the paper mentioning some conjectures and open questions in Section~\ref{sec:conj}.

\subsection{Permutations and consecutive patterns} We denote by $\S_n$ the set of permutations of $\{1,2,\dots,n\}$. If $\pi\in\S_n$,
we will write its one-line notation
as $\pi=[\pi(1),\pi(2),\dots,\pi(n)]$ (or $\pi=\pi(1)\pi(2)\dots\pi(n)$ if it creates no confusion). The use of square brackets is to distinguish it from
the cycle notation, where $\pi$ is written as a product of
cycles of the form $(i,\pi(i),\pi^2(i),\dots,\pi^{k-1}(i))$, with $\pi^k(i)=i$. For example, $\pi=[2,5,1,7,3,6,4]=(1,2,5,3)(4,7)(6)$.

Given a permutation $\pi=\pi(1)\pi(2)\dots\pi(n)$, let $D(\pi)$ denote the {\em descent set} of $\pi$, that is,
$D(\pi)=\{i\ :\ 1\le i\le n-1,\ \pi(i)>\pi(i+1)\}$. Let $\des(\pi)=|D(\pi)|$ be the number of descents.
The Eulerian polynomials are defined by  $$A_n(x)=\sum_{\pi\in\S_n} x^{\des(\pi)+1}.$$
Its coefficients are called the Eulerian numbers.
The descent set and the number of descents can be defined for any sequence of integers $a=a_1 a_2 \dots a_n$
by letting $D(a)=\{i\ :\ 1\le i\le n-1,\ a_i>a_{i+1}\}$.

\bs

Let $X$ be a totally ordered set, and let $x_1,\dots,x_n\in X$ with $x_1<x_2<\dots<x_n$. Any permutation of these values can be expressed as
$[x_{\pi(1)},x_{\pi(2)},\dots,x_{\pi(n)}]$, where $\pi\in\S_n$. We define its {\em reduction} to be $\rho([x_{\pi(1)},x_{\pi(2)},\dots,x_{\pi(n)}])=[\pi(1),\pi(2),\dots,\pi(n)]=\pi$.
Note that the reduction is just a relabeling of the entries with the numbers from $1$ to $n$, keeping the order relationships among them.
For example $\rho([4,7,1,6.2,\sqrt{2}])=[3,5,1,4,2]$. If the values $y_1,\dots,y_n$ are not all different, then $\rho([y_1,\dots,y_n])$ is not defined.

Given two permutations $\sigma\in\S_m$, $\pi\in\S_n$, with $m\ge n$, we say that $\sigma$ {\em contains} $\pi$ {\em as a consecutive pattern}
is there is some $i$ such that $\rho([\sigma(i),\sigma(i+1),\dots,\sigma(i+n-1)])=\pi$. Otherwise, we say that $\sigma$ {\em avoids}
$\pi$ {\em as a consecutive pattern}. The set of permutations in $\S_n$ that avoid $\pi$ as a consecutive pattern is denoted by $\A_n(\pi)$. We let
$\A(\pi)=\bigcup_{n\ge1} \A_n(\pi)$. Consecutive pattern containment was first studied in \cite{EliNoy}, where the sets $\A_n(\pi)$ are enumerated for certain
permutations $\pi$.

\subsection{Allowed and forbidden patterns in maps}\label{sec:forbdef}

Let $f$ be a map $f:X\rightarrow X$. Given $x\in X$ and $n\ge1$, we define
$$\Pat(x,f,n)=\rho([x,f(x),f^2(x),\dots,f^{n-1}(x)]),$$
provided that there is no pair $1\le i< j\le n$ such that
$f^{i-1}(x)=f^{j-1}(x)$. If there is such a pair, then $\Pat(x,f,n)$ is not defined. When it is defined, we have $\Pat(x,f,n)\in\S_n$.
If $\pi\in\S_n$ and there is some $x\in X$ such that $\Pat(x,f,n)=\pi$, we say that $\pi$ is {\em realized} by $f$ (at $x$), or that $\pi$ is an {\em allowed pattern} of $f$.
The set of all permutations realized by $f$ is denoted by $\Al(f)=\bigcup_{n\ge1} \Al_n(f)$, where $$\Al_n(f)=\{\Pat(x,f,n): x\in X\}\subseteq\S_n.$$
The remaining permutations are called {\em forbidden patterns}, and denoted by $\F(f)=\bigcup_{n\ge1} \F_n(f)$, where $\F_n(f)=\S_n\setminus\Al_n(f)$.

We are introducing some variations to the notation and terminology used in~\cite{AEK,Bandt}.
The main change is that our permutation $\pi=\Pat(x,f,n)$ is essentially the inverse of the permutation
of $\{0,1,\dots,n-1\}$ that the authors of~\cite{AEK} refer to as the {\em order pattern defined by $x$}. Our convention, aside from simplifying the notation,
will be more convenient from a combinatorial point of view. The advantage is that now the set $\Al(f)$  is closed under consecutive pattern containment,
in the standard sense used in the combinatorics literature.
Indeed, if $\sigma\in\Al(f)$ and $\sigma$ contains $\tau$ as a consecutive pattern, then $\tau\in\Al(f)$.
An equivalent statement is that if $\pi\in\F(f)$, then $\Al(f)\subseteq\A(\pi)$.
The set of minimal elements of $\F(f)$, i.e., those permutations in $\F(f)$ that avoid all other patterns in $\F(f)$, will be denoted $\B(f)$.

\ms

Let us consider now the case in which $X$ is a closed interval in $\mathbb{R}$, with the usual total order on real numbers.
An important incentive to study the set of forbidden patterns of a map comes from the following result, which is a consequence of \cite{Bandt}.
\begin{proposition}\label{prop:Bandt} If $I\subset \mathbb{R}$ is a closed interval and $f:I\rightarrow I$ is piecewise monotone, then $\F(f)\neq\emptyset$.
\end{proposition}
Recall that piecewise monotone means that there is a finite partition of $I$
into intervals such that $f$ is continuous and strictly monotone on each of
those intervals. It fact, it is shown in \cite{Bandt} that for such a map, the number of allowed patterns of $f$ grows at most exponentially, i.e.,
there is a constant $C$ such that $|\Al_n(f)|<C^n$ for $n$ large enough. The value of $C$ is related to the {\em topological entropy} of $f$ (see \cite{Bandt} for details).
Since the growth of the total number of permutations of length $n$ is super-exponential, the above proposition follows.

Proposition~\ref{prop:Bandt}, together with the above observation that $\Al(f)$ is closed under consecutive pattern containment,
provides an interesting connection between dynamical systems on one-dimensional interval maps and pattern avoiding permutations. An important application is that forbidden patterns can be used to distinguish
random from deterministic time series. Indeed, in a sequence $(x_1,x_2,x_3,\dots)$ where each $x_i$ has been chosen independently at random
from some continuous probability distribution, any pattern $\pi\in\S_n$ appears as $\pi=\rho([x_i,x_{i+1},\dots,x_{i+n-1}])$ for some $i$ with nonvanishing
probability, and this probability approaches one as the length of the sequence increases. On the other hand, if the sequence has been generated by defining
$x_{k+1}=f(x_k)$ for $k\ge1$, where $f:I\rightarrow I$ is a piecewise monotone map, then Proposition~\ref{prop:Bandt} guarantees that some patterns (in fact, most of them) will
never appear in the sequence. The practical aspect of these applications has been considered in~\cite{AZS}.

A new and interesting direction of research is to study more properties of the sets $\Al(f)$.
Some natural problems that arise are the following. \ben \renewcommand{\labelenumi}{\arabic{enumi}.}
\item Understand how $\Al(f)$ and $\B(f)$ depend on the map $f$.
\item\label{prob:enum} Describe and/or enumerate (exactly or asymptotically) $\Al(f)$ and $\B(f)$ for a particular~$f$.
\item Among the sets of patterns $\Sigma$ such that $\A_n(\Sigma)$ grows at most exponentially in $n$ (this is a necessary condition), characterize those for which there exists a map $f$ such that $\B(f)=\Sigma$.
\item Given a map $f$, determine the length of its smallest forbidden pattern.
\item\label{prob:periodic} Study how $\Al(f)$ is related to the set of permutations realized by periodic orbits of $f$, for which several results generalizing Sharkovskii's theorem are known (see \cite{Sarko,Bobok}).
\een
Most of this paper is devoted to solving problem~\ref{prob:enum} for a specific family of maps, that we describe next. Before that,
it is worth mentioning that problem~\ref{prob:periodic} is briefly considered
in~\cite{AEK}, where it is explained that for given $\pi\in\S_n$, the endpoints of the components of the set
$\{x\in I:\Pat(x,f,n)=\pi\}$ are points of period $n$ and (iterated) preimages of points of period less than $n$.

\subsection{One-sided shifts}

We will concentrate on the set of allowed patterns of certain maps called {\em one-sided shift maps}, or simply {\em shifts} for short. For a detailed definition of the associated dynamical system,
called the {\em one-sided shift space}, we refer the reader to~\cite{AEK}.

The totally ordered set $X$ considered above will now be the set $\W_N=\{0,1,\dots,N-1\}^\mathbb{N}$ of infinite words on $N$ symbols, equipped with
the lexicographic order.
Define the (one-sided) shift transformation \bce\bt{cccc} $\Sigma_N:$ & $\W_N$ & $\longrightarrow$ & $\W_N$ \\  & $w_1w_2w_3\dots$ & $\mapsto$ & $w_2w_3w_4\dots$.\et\ece
We will use $\Sigma$ instead of $\Sigma_N$ when it creates no confusion.

Given $w\in\W_N$, $n\ge1$, and $\pi\in\S_n$, we have from the above definition that $\Pat(w,\Sigma,n)=\pi$
if, for all indices $1\le i,j\le n$, $\Sigma^{i-1}(w)<\Sigma^{j-1}(w)$ if and only if $\pi(i)<\pi(j)$.
For example, \beq\label{eq:expat}\Pat(2102212210\dots,\Sigma,7)=[4,2,1,7,5,3,6],\eeq because the relative order of the successive shifts is
\bce\bt{rll}
$2102212210\dots$ & \quad & $4$ \\
$102212210\dots$ & \quad & $2$ \\
$02212210\dots$ & \quad & $1$ \\
$2212210\dots$ & \quad & $7$ \\
$212210\dots$ & \quad & $5$ \\
$12210\dots$ & \quad & $3$ \\
$2210\dots$ & \quad & $6$, \\ \et\ece
regardless of the entries in place of the dots.
The case $N=1$ is trivial, since the
only allowed pattern of $\Sigma_1$ is the permutation of length 1. In the rest of the paper, we will assume that $N\ge2$.

%One of our goals is to study the set $\F(\Sigma_N)$, the set of {\em forbidden patterns} of $\Sigma_N$.
%We will also be interested in the set $\R(\Sigma_N)$ of minimal forbidden patterns with respect to pattern containment.

If $x\in\{0,1,\dots,N-1\}$, we will use the notation $x^\infty=xxx\dots$. If $w\in\W_N$, then $w_n$ denotes the $n$-th letter of $w$, and we write $w=w_1w_2w_3\dots$.
We will also write $w_{[k,\ell]}=w_kw_{k+1}\dots w_{\ell}$, and $w_{[k,\infty)}=w_kw_{k+1}\dots$.
Note that $w_{[k,\infty)}=\Sigma^{k-1}(w)$.

\begin{figure}[hbt] \centering
\epsfig{file=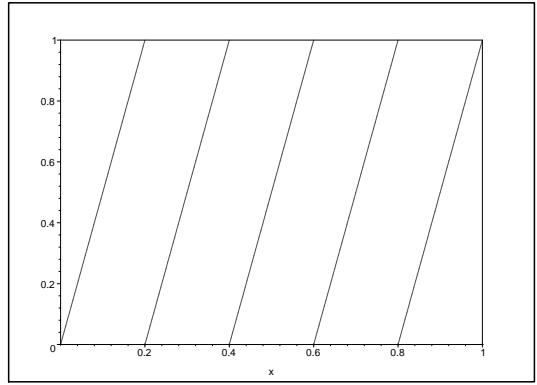,height=7cm,angle=-90} \caption{The sawtooth map $x \mapsto Nx\mod 1$ for $N=5$.}\label{fig:sawtooth}
\end{figure}

Let $\Upsilon\subset\W_N$ be the set of all words of the form $u(N{-}1)^\infty$ except the word $(N{-}1)^\infty$. As shown in~\cite{AEK}, $\W_N\setminus\Upsilon$ is closed under one-sided shifts, and the map
\bce\bt{cccc} $\varphi:$ & $\W_N\setminus\Upsilon$ & $\longrightarrow$ & $[0,1]$ \\  & $w_1w_2w_3\dots$ & $\mapsto$ & $\sum_{i\ge0} w_i N^{-(i+1)}$\et\ece
is an order-preserving isomorphism (here $[0,1]$ denotes the closed unit interval). The map $$\varphi\circ\Sigma_N\circ\varphi^{-1}:[0,1]\rightarrow[0,1]$$ is the so-called {\em sawtooth map} $x \mapsto Nx\mod 1$, which
is piecewise linear (see Figure~\ref{fig:sawtooth}). By Proposition~\ref{prop:Bandt}, this map has forbidden patterns.
It is not hard to check (see~\cite{AEK}) that
$\Sigma_N$ has the same set of forbidden patterns as this map, and in particular $\F(\Sigma_N)\neq\emptyset$.
Forbidden patterns of shift systems were first studied in~\cite{AEK}, where the main result is that the smallest forbidden patterns of $\Sigma_N$
have length $N+2$.

\begin{proposition}[\cite{AEK}]\label{th:aek_minshift}Let $N\ge2$. We have that
\renewcommand{\labelenumi}{(\alph{enumi})}
\ben \item $\F_n(\Sigma_N)=\emptyset$ for every $n\le N+1$,
\item $\F_n(\Sigma_N)\neq\emptyset$ for every $n\ge N+2$. %In fact, $\R_n(\Sigma_N)\neq\emptyset$ in this case.
\een
\end{proposition}

\begin{example} One can check that the smallest forbidden patterns of $\Sigma_4$ are
$$615243, 324156, 342516, 162534, 453621, 435261.$$
\end{example}

Recall that a word $w\in\{0,1,\dots,N-1\}^k$ is {\em primitive} if it cannot be written as a power of any proper subword, i.e.,
it is not of the form $w=u^m$ for any $m>1$, where the exponent indicates concatenation of $u$ with itself $m$ times.
Let $\psi_N(k)$ denote the number of primitive words of length $k$ over an $N$-letter alphabet. It is well known that $$\psi_N(k)=\sum_{d|k}\mu(d)N^{k/d},$$
where $\mu$ denotes the M\"obius function.

\section{The number of symbols needed to realize a pattern}
\label{sec:numsymb}

Given a permutation $\pi\in\S_n$, let $N(\pi)$ be the smallest number $N$ such that $\pi\in\Al(\Sigma_N)$. The value of $N(\pi)$ indicates
what is the minimum number of symbols needed in the alphabet in order for $\pi$ to be realized by a shift. For example, if $\pi=[4,2,1,7,5,3,6]$, then $N(\pi)\le3$ because of~(\ref{eq:expat}),
and it is not hard to see that $N(\pi)=3$. The main result in this section is a formula for $N(\pi)$.

\begin{theorem}\label{th:Npi} Let $n\ge2$. For any $\pi\in\S_n$, $N(\pi)$ is given by
\beq N(\pi)=1+|A(\pi)|+\Delta(\pi),\label{eq:Npi} \eeq
where $$A(\pi)=\{a\,:\,1\le a\le n-1 \mbox{ such that if } i=\pi^{-1}(a),\,j=\pi^{-1}(a+1), \mbox{ then } i,j<n \mbox{ and } \pi(i+1)>\pi(j+1)\},$$
and
$\Delta(\pi)=0$ except in the following three cases, in which $\Delta(\pi)=1$:
\ben\renewcommand{\labelenumi}{(\Roman{enumi})}
\item $\pi(n)\notin\{1,n\}$, and if $i=\pi^{-1}(\pi(n)-1)$, $j=\pi^{-1}(\pi(n)+1)$, then $\pi(i+1)>\pi(j+1)$;
\item $\pi(n)=1$ and $\pi(n-1)=2$; or
\item $\pi(n)=n$ and $\pi(n-1)=n-1$.
\een
\end{theorem}

Note that $A(\pi)$ is the set of entries $a$ in the one-line notation of $\pi$ such that the entry following $a+1$ is smaller than the entry following $a$.
For example, if $\pi=[4, 3, 6, 1, 5, 2]$, then $A(\pi)=\{3,4,5\}$, so Theorem~\ref{th:Npi} says that $N(\pi)=1+3+0=4$.
The following lemmas will be useful in the proof.

\begin{lemma}\label{lem:sym}
Suppose that $\Pat(w,\Sigma,n)=\pi$, where $w\in\W_N$ and $\pi\in\S_n$. Define $w'=w'_1w'_2\dots$ by $w'_i=N-1-w_i$ for $i\ge1$, and $\pi'=[\pi'(1),\dots,\pi'(n)]$ by $\pi'(j)=n+1-j$ for $1\le j\le n$.
Then, $\Pat(w',\Sigma,n)=\pi'$.
\end{lemma}

\begin{proof}
It is clear from the definition of $w'$ that for any $i,j\ge 1$, $\Sigma^{i-1}(w)<\Sigma^{j-1}(w)$ if and only if $\Sigma^{i-1}(w')>\Sigma^{j-1}(w')$. Thus, $\Pat(w',\Sigma,n)$ is the only permutation $\pi'\in\S_n$ that satisfies that
$\pi(i)<\pi(j)$ if and only if $\pi'(i)>\pi'(j)$, for all $1\le i,j\le n$.
\end{proof}

\begin{lemma}\label{lemma:1}
Suppose that $\Pat(w,\Sigma,n)=\pi$. If $1\le i,j<n$, $\pi(i)<\pi(j)$, and $\pi(i+1)>\pi(j+1)$, then $w_i<w_j$.
\end{lemma}

\begin{proof}
Clearly $\pi(i+1)>\pi(j+1)$ implies that $w_{[i+1,\infty)}>w_{[j+1,\infty)}$ lexicographically. Similarly, $\pi(i)<\pi(j)$ implies that
$w_{[i,\infty)}=w_iw_{[i+1,\infty)}<w_jw_{[j+1,\infty)}=w_{[j,\infty)}$. This can only happen if $w_i<w_j$.
\end{proof}

\begin{lemma}\label{lem:primitive}
Suppose that $\Pat(w,\Sigma,n)=\pi$. If $1\le i<k\le n$ are such that $|\pi(i)-\pi(k)|=1$, then the word $w_{[i,k-1]}$ is primitive.
\end{lemma}

\begin{proof}
If $w_{[i,k-1]}$ is not primitive, we can write $w_{[i,k-1]}=q^r$ for some $r>1$. Let $v=w_{[k,\infty)}$, and let $j=i+|q|$.
We have that $w_{[i,\infty)}=q^r v$, $w_{[j,\infty)}=q^{r-1} v$, and $w_{[k,\infty)}=v$.

If $\pi(i)<\pi(j)$, then $q^r v<q^{r-1} v$, which, by canceling the prefix $q^{r-l}$, implies that $q^l v<q^{l-1} v$ for any $1\le l \le r$.
By iteration of this we have that $q^{r-1} v< v$, and thus $\pi(i)<\pi(j)<\pi(k)$, which is a contradiction with $|\pi(i)-\pi(k)|=1$.
If $\pi(i)>\pi(j)$, then the same argument with the inequalities reversed shows that $q^{r-1} v> v$, so $\pi(i)>\pi(j)>\pi(k)$, which is a contradiction as well.
\end{proof}

We will prove Theorem~\ref{th:Npi} in two parts. In Subsection~\ref{sec:lowerbound} we prove that $1+|A(\pi)|+\Delta(\pi)$ is a lower bound
on $N(\pi)$, and in Subsection~\ref{sec:upperbound} we prove that it is an upper bound as well.

\subsection{Proof of $N(\pi)\ge 1+|A(\pi)|+\Delta(\pi)$.}\label{sec:lowerbound}
To show that the right hand side of~(\ref{eq:Npi}) is a lower bound on $N(\pi)$,
we will prove that if $w\in\W_N$ is such that $\Pat(w,\Sigma,n)=\pi$, then necessarily
$N\ge 1+|A(\pi)|+\Delta(\pi)$. This fact is a consequence of the following lemma.

\begin{lemma}\label{lem:required}
Suppose that $\Pat(w,\Sigma,n)=\pi$, and let $b=\pi(n)$. The entries of $w$ satisfy
\beq \label{eq:increasing} w_{\pi^{-1}(1)}\le w_{\pi^{-1}(2)} \le \dots \le w_{\pi^{-1}(n)}, \eeq
with strict inequalities $w_{\pi^{-1}(a)}< w_{\pi^{-1}(a+1)}$ for each $a\in A(\pi)$. Additionally, if $\Delta(\pi)=1$, then in each of the three cases from Theorem~\ref{th:Npi} we have, respectively, that
\ben\renewcommand{\labelenumi}{(\Roman{enumi})}
\item one of the inequalities $w_{\pi^{-1}(b-1)}\le w_{n}\le w_{\pi^{-1}(b+1)}$ is strict;
\item $\dots\le  w_{n+2} \le w_{n+1} \le w_n \le w_{n-1}$ and one of these inequalities is strict;
\item $w_{n-1}\le w_n\le w_{n+1}\le w_{n+2}\le \cdots$ and one of these inequalities is strict.
\een
In all cases, the entries of $w$ must satisfy $|A(\pi)|+\Delta(\pi)$ strict inequalities.

\end{lemma}

\begin{proof}
The condition $\Pat(w,\Sigma,n)=\pi$ is equivalent to
\beq \label{eq:condition} w_{[\pi^{-1}(1),\infty)}< w_{[\pi^{-1}(2),\infty)} < \dots < w_{[\pi^{-1}(n),\infty)}, \eeq
which clearly implies equation~(\ref{eq:increasing}).
If we remove the term $w_{n}$ from~(\ref{eq:increasing}), we get
\beq\label{eq:increasing_n} \begin{cases}
\mbox{(a)} \quad w_{\pi^{-1}(1)}\le w_{\pi^{-1}(2)} \le \dots \le w_{\pi^{-1}(b-1)}\le w_{\pi^{-1}(b+1)} \le \dots \le w_{\pi^{-1}(n)} & \mbox{ if } b\notin\{1,n\}, \\
\mbox{(b)} \quad w_{\pi^{-1}(2)} \le w_{\pi^{-1}(3)} \le\dots \le w_{\pi^{-1}(n)} & \mbox{ if } b=1, \\
\mbox{(c)} \quad w_{\pi^{-1}(1)} \le w_{\pi^{-1}(2)} \le\dots \le w_{\pi^{-1}(n-1)} & \mbox{ if } b=n.
\end{cases}\eeq

For every $a\in A(\pi)$, the inequality $w_{\pi^{-1}(a)}< w_{\pi^{-1}(a+1)}$ in (\ref{eq:increasing_n}) has to be strict, by Lemma~\ref{lemma:1} with
$i=\pi^{-1}(a)$ and $j=\pi^{-1}(a+1)$.
Let us now see that in the three cases when $\Delta(\pi)=1$, an additional strict inequality must be satisfied.

Consider first case (I). Let $i=\pi^{-1}(b-1)$ and $j=\pi^{-1}(b+1)$.
Since $\pi(i+1)>\pi(j+1)$, Lemma~\ref{lemma:1} implies that $w_i<w_j$, so the inequality
$w_{\pi^{-1}(b-1)}< w_{\pi^{-1}(b+1)}$ (equivalently, $w_i<w_j$) in (\ref{eq:increasing_n}a) has to be strict.
In case (II), the leftmost inequality in~(\ref{eq:condition}) is $w_{[n,\infty)} < w_{[n-1,\infty)}$. For this to hold, we need
$\dots\le  w_{n+2} \le w_{n+1} \le w_n \le w_{n-1}$ and at least one of these inequalities must be strict.
Similarly, in case (III), the rightmost inequality in~(\ref{eq:condition}) is $w_{[n-1,\infty)} < w_{[n,\infty)}$. This forces $w_{n-1}\le w_n\le w_{n+1}\le w_{n+2}\le\cdots$
with at least one strict inequality.
\end{proof}

We will refer to the $|A(\pi)|+\Delta(\pi)$ strict inequalities in Lemma~\ref{lem:required} as the {\em required strict inequalities}. Combined with the weak inequalities from the lemma,
they force the number of symbols used in $w$ to be at least $1+|A(\pi)|+\Delta(\pi)$. Examples~2 and~3 illustrate how this lemma is used.

\subsection{Proof of $N(\pi)\le 1+|A(\pi)|+\Delta(\pi)$.}\label{sec:upperbound}

We will show how for any given $\pi\in\S_n$ one can construct a word $w\in\W_{N}$ with
$\Pat(w,\Sigma,n)=\pi$, where $N=1+|A(\pi)|+\Delta(\pi)$. We need $w$ to satisfy condition~(\ref{eq:condition}). Again, let $b=\pi(n)$.

The first important observation is that, if we can only use $N$ different symbols, then the
$|A(\pi)|+\Delta(\pi)=N-1$ required strict inequalities from Lemma~\ref{lem:required} determine the values of the entries $w_1w_2\dots w_{n-1}$. This fact is
restated as Corollary~\ref{cor:determined}. Consequently, we are forced to assign values to these entries as follows:

\ben \renewcommand{\labelenumi}{(\alph{enumi})}
\item If $b\notin\{1,n\}$, assign values to the variables in~(\ref{eq:increasing_n}a) from left to right, starting with $w_{\pi^{-1}(1)}=0$ and increasing the value by 1 at each required strict inequality.
\item If $b=1$, assign values to the variables in~(\ref{eq:increasing_n}b) from left to right, starting with $w_{\pi^{-1}(2)}=0$ if $\pi(n-1)\neq 2$, or with $w_{\pi^{-1}(2)}=1$ if $\pi(n-1)=2$
(this is needed in order for condition (II) in Lemma~\ref{lem:required} to hold),
and increasing the value by 1 at each required strict inequality.
\item If $b=n$, assign values to the variables in~(\ref{eq:increasing_n}c) from left to right, starting with $w_{\pi^{-1}(1)}=0$ and increasing the value by 1
at each required strict inequality. (Note that when $\Delta(\pi)=1$, the last assigned value is $w_{\pi^{-1}(n-1)}=w_{n-1}=|A(\pi)|=N-2$.)
\een

It remains to assign the values to $w_m$ for $m\ge n$. Before we do this, let us prove some facts about the entries $w_1\dots w_{n-1}$. In the
following three lemmas, $\pi$ is any permutation in $\S_n$ with $N(\pi)=N$ and $w_1\dots w_{n-1}$ are the values in $\{0,1,\dots,N-1\}$ assigned above in order to satisfy the required strict inequalities.

\begin{lemma}\label{lem:nonzero}
Let $i<n$. If $\pi(i)>\pi(i+1)$, then $w_i\ge1$.
\end{lemma}

\begin{proof}
In the case that $\pi(n)=1$ and $\pi(n-1)=2$, our construction makes $w_j\ge 1$ for $1\le j\le n-1$, so this case is trivial.

In any other case, we will show that in equation~(\ref{eq:increasing_n})
there is some required strict inequality to the left of $w_{i}$.
Let $c=\pi(i)$. Consider the sequence $1,2,\dots,c$ and remove the entry $\pi(n)$ if it is in the sequence.
Let $a_1,\dots,a_s$ be the resulting sequence (note that $s=c-1$ or $s=c$ depending on whether $\pi(n)\le c$ or not).
For $1\le \ell\le s$, let $b_\ell$ be the entry following $a_\ell$ in the one-line notation of $\pi$.
We have that $a_1<b_1$ and $a_{s}=c>b_{s}=\pi(i+1)$. Since $a_{\ell+1}-a_\ell\le 2$ for every $\ell$ and all the $b_\ell$'s are different,
there has to be some $r$ such that $b_r>b_{r+1}$.
Now the inequality $w_{\pi^{-1}(a_r)}<w_{\pi^{-1}(a_{r+1})}$ must be one of the required strict inequalities,
because either it is of the form $w_{\pi^{-1}(a)}< w_{\pi^{-1}(a+1)}$
with $a\in A(\pi)$, or it is $w_h< w_j$ with $h=\pi^{-1}(b-1)$, $j=\pi^{-1}(b+1)$, and $\pi(h+1)>\pi(j+1)$.
\end{proof}

\begin{lemma}\label{lem:nonmax}
Let $i<n$. If $\pi(i)<\pi(i+1)$, then $w_i\le N-2$.
\end{lemma}

\begin{proof}
In the case that $\pi(n)=n$ and $\pi(n-1)=n-1$, it follows from our construction that $w_j\le |A(\pi)|=N-2$ for $1\le j\le n-1$, so this case is trivial.

In any other case, we will show that in equation~(\ref{eq:increasing_n})
there is some required strict inequality to the right of $w_{i}$.
Let $c=\pi(i)$. Consider the sequence $c,c+1,\dots,n$ and remove the entry $\pi(n)$ if it is in the sequence.
Let $a_1,\dots,a_s$ be the resulting sequence.
For $1\le \ell\le s$, let $b_\ell$ be the entry following $a_\ell$ in the one-line notation of $\pi$.
We have that $a_1=c<b_1=\pi(i+1)$ and $a_{s}>b_{s}$. Since $a_{\ell+1}-a_\ell\le 2$ for every $\ell$ and all the $b_\ell$'s are different,
there has to be some $r$ such that $b_r>b_{r+1}$.
Now the inequality $w_{\pi^{-1}(a_r)}<w_{\pi^{-1}(a_{r+1})}$ must be one of the required strict inequalities,
as in the proof of Lemma~\ref{lem:nonzero}.
\end{proof}

\begin{lemma}\label{lem:order}
If $1\le i,j<n$ are such that $\pi(i)<\pi(j)$ and $\pi(i+1)>\pi(j+1)$, then $w_i<w_j$.
\end{lemma}

\begin{proof}
We will show that in equation~(\ref{eq:increasing_n})
there is some required strict inequality between $w_i$ and $w_j$.
Let $a_1,\dots,a_{s}$ be the ordered sequence of integers between $\pi(i)$ and $\pi(j)$, both included,
after $\pi(n)$ has been removed from it (if it appeared).
For $1\le \ell\le s$, let $b_\ell$ be the entry following $a_\ell$ in the one-line notation of $\pi$.
The fact that $\pi(i+1)=b_1>b_\ell=\pi(j+1)$ implies that there has to be some $1\le r<s$ such that $b_r>b_{r+1}$.
Now the inequality $w_{\pi^{-1}(a_r)}<w_{\pi^{-1}(a_{r+1})}$ has to be strict,
as in the proof of Lemma~\ref{lem:nonzero}.
\end{proof}

\ms

Once the values $w_1\dots w_{n-1}$ have been determined, there are several ways to assign values to $w_m$ for $m\ge n$. Here we give a few different possibilities. Each one has some restriction on the value of $b=\pi(n)$.

\ben\renewcommand{\labelenumi}{\Alph{enumi}.}
\item Assume that $b\neq n$. Let $k=\pi^{-1}(b+1)$. Let $u=w_1w_2\dots w_{k-1}$ and $p=w_kw_{k+1}\dots w_{n-1}$.
Let $m$ be any integer satisfying $m \ge 1+\frac{n-2}{n-k}$ (for definiteness, we can pick $m=n-1$).
Let $w_A(\pi)=up^{m}0^\infty$.

\item Assume that $b\neq 1$. Let $k=\pi^{-1}(b-1)$. Let $u=w_1w_2\dots w_{k-1}$ and $p=w_kw_{k+1}\dots w_{n-1}$.
Again, let $m$ be such that $m \ge 1+\frac{n-2}{n-k}$ (for definiteness, we can pick $m=n-1$). Let $w_B(\pi)=up^{m}(N{-}1)^\infty$.

\item Assume that $b=1$. Let $w_C(\pi)=w_1w_2\dots w_{n-1}0^\infty$.
\item Assume that $b=n$. Let $w_D(\pi)=w_1w_2\dots w_{n-1}(N{-}1)^\infty$.

\item Assume that $1<b<n$, and that $w_{\pi^{-1}(b-1)}< w_{\pi^{-1}(b+1)}$ (this happens when $\Delta(\pi)=1$). Let $c=w_{\pi^{-1}(b-1)}$.
Define $w_E(\pi)=w_1w_2\dots w_{n-1}c(N{-}1)^\infty$ and $w_F(\pi)=w_1w_2\dots w_{n-1}(c{+}1)0^\infty$.
\een

Clearly, $w=w_\ast(\pi)$ uses $N$ different symbols in each one of the above cases. It remains to prove that $\Pat(w,\Sigma,n)=\pi$, which is equivalent
to showing that $w$ satisfies condition~(\ref{eq:condition}). Let us now prove that this is the case for $w=w_A(\pi)$, when $b\neq n$.

In the following three lemmas and in Proposition~\ref{prop:wworks}, $\pi$ is any permutation in $\S_n$ with $\pi(n)\neq n$, and $w=w_A(\pi)$. Also, $k$, $u$, $p$ and $m$ are as defined in case A above.

\begin{lemma}\label{lem:non0}
The word $p=w_kw_{k+1}\dots w_{n-1}$ has some nonzero entry.
\end{lemma}

\begin{proof}
Since $\pi(k)=\pi(n)+1>\pi(n)$, there has to be a descent in the sequence $\pi(k),\pi(k+1),\dots,\pi(n-1),\pi(n)$. Let $k\le i\le n-1$ be such that $\pi(i)>\pi(i+1)$.
By Lemma~\ref{lem:nonzero}, $w_i\ge1$, so we are done.
\end{proof}

\begin{lemma}\label{lem:pprimitive}
The word $p=w_kw_{k+1}\dots w_{n-1}$ is primitive.
\end{lemma}

\begin{proof}
Assume for contradiction that $p=q^r$ for some $r>1$. Let $s=|q|$. Then, $k\le i<n-s$, we have that $w_i=w_{i+s}$.
There are now two possibilities depending on the relationship between $\pi(k)$ and $\pi(k+s)$.

If $\pi(k)<\pi(k+s)$, then by Lemma~\ref{lem:order}, $\pi(k+1)<\pi(k+s+1)$. By iterating this argument, we see that $\pi(i)<\pi(i+s)$ for $k\le i\le n-s$. In particular,
$\pi(k)<\pi(k+s)<\pi(k+2s)<\dots<\pi(k+rs)=\pi(n)$, which contradicts the choice of $k$ as the index such that $\pi(k)=\pi(n)+1$.

If $\pi(k)>\pi(k+s)$, then similarly by Lemma~\ref{lem:order} we have that $\pi(i)>\pi(i+s)$ for $k\le i\le n-s$. In particular,
$\pi(k)>\pi(k+s)>\dots>\pi(k+rs)=\pi(n)$. Since $r>1$, this implies that $\pi(k)-\pi(n)\ge2$, which contradicts again the choice of $k$.
\end{proof}

\begin{lemma}\label{lem:nk} We have that $w_{[n,\infty)}<w_{[k,\infty)}$. Moreover, there is no $1\le s\le n$ such that $w_{[n,\infty)}<w_{[s,\infty)}<w_{[k,\infty)}$.
\end{lemma}

\begin{proof}
For the first part, note that $w_{[n,\infty)}=p^{m-1}0^\infty<w_{[k,\infty)}=p^{m}0^\infty$ because $p$ has some nonzero entry, by Lemma~\ref{lem:non0}.
 For the second part, assume there is an $s$ such that $p^{m-1}0^\infty<w_{[s,\infty)}<p^{m}0^\infty$. Then $w_{[s,\infty)}=p^{m-1}v$, with $0^\infty<v<p0^\infty$.
Since $p$ is primitive (by Lemma~\ref{lem:pprimitive}), the initial $p$ in $w_{[s,\infty)}$ cannot overlap simultaneously with the two occurrences of $p$ at the beginning of $w_{[k,\infty)}$. Given that $s\neq k,n$, this only
leaves the possibility that $s<k$. If some of the occurrences of $p$ in $w_{[s,\infty)}=p^{m-1}v$ and in
$w_{[k,\infty)}=p^{m}0^\infty$ coincide, then the condition $v<p0^\infty$ does not hold. The only remaining possibility is that the $m-1$ initial occurrences of $p$ in $w_{[s,\infty)}$
are entirely contained in $u=w_1w_2\dots w_{k-1}$.
However, the choice of $m$ satisfying $m\ge 1+\frac{n-2}{n-k}$ was made so that $|p^{m-1}|=(n-k)(m-1)\ge n-2$. This only leaves the case that $s=1$ and $u=p^{m-1}$,
but then $v=w_{[k,\infty)}=p^m0^\infty$, so the condition $v<p0^\infty$ does not hold.
\end{proof}

The next proposition proves that $\Pat(w_A(\pi),\Sigma,n)=\pi$.

\begin{proposition}\label{prop:wworks}
If $1\le i,j\le n$ are such that $\pi(i)<\pi(j)$, then $w_{[i,\infty)}<w_{[j,\infty)}$.
\end{proposition}

\begin{proof}
Let $S(i,j)$ be the statement ``$\pi(i)<\pi(j)$ implies $w_{[i,\infty)}<w_{[j,\infty)}$". We want to prove $S(i,j)$ for all $1\le i,j\le n$ with $i\neq j$.
We consider three cases.

\bit \item {\it Case $i=n$}. Assume that $\pi(n)<\pi(j)$. By Lemma~\ref{lem:nk} we know that $w_{[n,\infty)}<w_{[k,\infty)}$. If $j=k$ we are done. If $j\ne k$, then $\pi(n)<\pi(j)$
implies that $\pi(k)<\pi(j)$, since $\pi(k)=\pi(n)+1$. So, if $S(k,j)$ holds, then $S(n,j)$ must hold as well. We have reduced $S(n,j)$ to $S(k,j)$. Equivalently, $\neg S(n,j)\Rightarrow\neg S(k,j)$, where $\neg$ denotes negation.

\item {\it Case $j=n$}. Assume that $\pi(i)<\pi(n)$, so in particular $i\ne k$. By the second part of Lemma~\ref{lem:nk}, in order to prove that $w_{[i,\infty)}<w_{[n,\infty)}$
it is enough to show that $w_{[i,\infty)}<w_{[k,\infty)}$. Also, $\pi(i)<\pi(n)$ implies that $\pi(i)<\pi(k)$. We have reduced $S(i,n)$ to $S(i,k)$.

\item {\it Case $i,j<n$}. Assume that $\pi(i)<\pi(j)$. If $w_i<w_j$, then $w_{[i,\infty)}<w_{[j,\infty)}$ and we are done. If $w_i=w_j$, then we know by Lemma~\ref{lem:order} that
$\pi(i+1)<\pi(j+1)$. If we can show that $w_{[i+1,\infty)}<w_{[j+1,\infty)}$, then $w_{[i,\infty)}=w_iw_{[i+1,\infty)}<w_jw_{[j+1,\infty)}=w_{[j,\infty)}$.
So, we have reduced $S(i,j)$ to $S(i+1,j+1)$.
\eit

The above three cases prove that for all $1\le i,j\le n-1$, $\neg S(i,j)\Rightarrow\neg S(g(i),g(j))$, where $g$ is defined for $1\le i\le n-1$ as
$$g(i)=\begin{cases} i+1 & \mbox{if } i<n-1, \\ k & \mbox{if } i=n-1.\end{cases}$$
Suppose now that $S(i,j)$ does not hold for some $i,j$. Using the first two cases above, we can assume that $1\le i,j\le n-1$. Then, $S(g^\ell(i),g^\ell(j))$ fails for every $\ell\ge 1$.
Let $r$ be the index for which $\pi(r)$ is maximum among those $r$ with $k\le r\le n-1$. Let $\ell$ be such that $g^\ell(i)=r$
and $k\le g^\ell(j)\le n-1$ (such an $\ell$ always exists).
Then $S(g^\ell(i),g^\ell(j))$ must hold because $\pi(g^\ell(i))\ge \pi(g^\ell(j))$. This is a contradiction, so the proposition is proved.
\end{proof}

If $b\neq 1$, proving that $w=w_B(\pi)$ satisfies $\Pat(w,\Sigma,n)=\pi$ is analogous to our argument for $w=w_A(\pi)$. Instead of Lemma~\ref{lem:nonzero} we would use
Lemma~\ref{lem:nonmax}, and the analogue of Lemma~\ref{lem:non0} is the fact that $p=w_kw_{k+1}\dots w_{n-1}$ has some entry strictly smaller than $N-1$. Lemma~\ref{lem:nk} would be replaced with the fact that
$w_{[k,\infty)}<w_{[n,\infty)}$ and there is no $1\le s\le n$ such that $w_{[k,\infty)}<w_{[s,\infty)}<w_{[n,\infty)}$. With these variations, Proposition~\ref{prop:wworks} can be proved similarly for $w=w_B(\pi)$.

We can complete the proof of the upper bound on $N(\pi)$ as follows. Let $\pi\in\S_n$ be given, and let $N=1+|A(\pi)|+\Delta(\pi)$. If $\pi(n-1)>\pi(n)$, let $w=w_A(\pi)$.
If $\pi(n-1)<\pi(n)$, let $w=w_B(\pi)$. Since $\Pat(w,\Sigma,n)=\pi$ and $w\in\W_N$, the theorem is proved.

\bs

\begin{example}
Let $\pi=[4, 3, 6, 1, 5, 2]$. By Theorem~\ref{th:Npi}, $N(\pi)=4$.
If $\Pat(w,\Sigma,n)=\pi$, then Lemma~\ref{lem:required} implies that
$w_4\le w_6\le w_2 < w_1 < w_5 < w_3$, and
there are no more required strict inequalities. We assign $w_4=w_2=0$, $w_1=1$, $w_5=2$, $w_3=3$.
Since $\pi(5)>\pi(6)$ and $b=\pi(6)=2$, we can take $w=w_A(\pi)$ (with $m=2$), so $k=\pi^{-1}(3)=2$, $u=w_1=1$, and $p=w_2w_3w_4w_5=0302$. We get $w=up^20^\infty=1030203020^\infty$.
\end{example}

The following consequence of the proof of Theorem~\ref{th:Npi} will be used in Section~\ref{sec:binary}.

\begin{corollary}\label{cor:determined}
Let $\pi\in\S_n$, $N=N(\pi)$, and let $w\in\W_N$ be such that $\Pat(w,\Sigma,n)=\pi$.
Then the entries $w_1w_2\dots w_{n-1}$ are uniquely determined by $\pi$.
\end{corollary}

Note that, however, with the conditions of Corollary~\ref{cor:determined},
$w_n$ is not always determined. In the case that $\pi(n)\notin\{1,n\}$ and $\Delta(\pi)=1$, we have two choices for $w_n$. In general,
there is a lot of flexibility in the choice of $w_m$ for $m\ge n$. The choices $w=w_A(\pi)$ and $w=w_B(\pi)$ in the proof of Theorem~\ref{th:Npi}
were made to simplify the proof of Proposition~\ref{prop:wworks} for all cases at once. It is not hard to modify the proof to see that $w_C(\pi)$, $w_D(\pi)$, $w_E(\pi)$ and $w_F(\pi)$
work as well, in the cases where they are defined. Let us see some examples of different ways to find a word $w\in\W_N$ such that $\Pat(w,\Sigma,n)=\pi$, for given $\pi\in\S_n$ with $N(\pi)=N$.

\begin{example} \renewcommand{\labelenumi}{(\roman{enumi})}\ben \item Let $\pi=[8,9,3,1,4,6,2,7,5]$. If $\Pat(w,\Sigma,n)=\pi$, then
by Lemma~\ref{lem:required} we must have $w_4\le w_7 < w_3 \le w_5 \le w_9 \le w_6 \le w_8 \le w_1 < w_2$, with a strict inequality $w_5<w_6$
caused by the fact that $\Delta(\pi)=1$. We have $N(\pi)=4 = 1+|\{2,8\}|+1$.
If $w\in\W_4$, the entries $w_4=w_7=0$, $w_3=w_5=1$, $w_6=w_8=w_1=2$, $w_2=3$ are forced. For the remaining entries we could take for example $w_{[9,\infty)}=13^\infty$, obtaining $w=w_E(\pi)=231012021 3^\infty$.
\item If $\pi=[3,5,2,4,1]$, then  Lemma~\ref{lem:required} tells us that $w_5\le w_3\le w_1< w_4\le w_2$ and there are no more required strict inequalities. We can take $w=w_C(\pi)=01010 0^\infty$.
\item If $\pi=[2,3,1,4,5,6]$, then by  Lemma~\ref{lem:required} we have that $w_3< w_1< w_2\le w_4\le w_5\le w_6$, and there must be some strict inequality in $w_5\le w_6\le w_7 \le \cdots$.
We can take $w=w_D(\pi)=120223 3^\infty$.
\een
\end{example}

Note that in terms of $\pi^{-1}$, the set $A(\pi)$ can alternatively be described as
$$A(\pi)=\{i\,:\, \pi^{-1}(i)+1 \mbox{ appears to the right of }\pi^{-1}(i+1)+1\mbox{ in the one-line notation for }\pi^{-1}\}.$$
For instance, if $\pi$ is the permutation in part (i) of the above example, we have $\pi^{-1}=[4,7,3,5,9,6,8,1,2]$, so $A(\pi)=\{2,8\}$.

\section{An equivalent characterization} %Permutations requiring the most number of symbols}
\label{sec:mostsymb}

We start this section by giving an expression for $N(\pi)$ that is sometimes more convenient to work with than the one in Theorem~\ref{th:Npi}.
We denote by $\mathcal{C}_n$ the set of permutations in $\S_n$ whose cycle decomposition consists of a unique cycle of length $n$.
Let $\T_n$ be the set of permutations $\pi\in\mathcal{C}_n$ with one distinguished entry $\pi(i)$, for some $1\le i\le n$. We call the elements of $\T_n$ {\em marked cycles}. We will use the symbol $\star$ to
denote the distinguished entry, both in one-line and in cycle notation. Note that it is not necessary to keep track of its value, since it is determined once we know all the remaining entries.
For example, $\T_3=\{[\star,3,1],[2,\star,1],[2,3,\star],[\star,1,2],[3,\star,2],[3,1,\star]\}$. Clearly, $|\T_n|=(n-1)!\cdot n=n!$, since there are $n$ ways to choose the distinguished entry.

Define a map \bea\nn \theta:\ \S_n&\rightarrow&\T_n \\ \nn \pi & \mapsto & \hp \eea as follows.
Given $\pi=[\pi(1),\pi(2),\dots,\pi(n)]$, let $\hp$ be the permutation
whose cycle decomposition is $(\pi(1),\pi(2),\dots,\pi(n))$, with the entry $\pi(1)$ distinguished, that is, $\hp=(\star,\pi(2),\dots,\pi(n))$. For example, if $\pi=[8,9,2,3,6,4,1,5,7]$, then $\hp=(\star,9,2,3,6,4,1,5,7)$, or in one-line notation, $\hp=[5,3,6,1,7,4,\star,9,2]$.
 Equivalently, for each $1\le i\le n$, we could define $\hp(i)$ to be the entry immediately to the right of $i$ in the one-line notation of $\pi$ if $i\neq\pi(n)$, and $\hp(i)=\star$ if $i=\pi(n)$.

The map $\theta$ is a bijection between $\S_n$ and $\T_n$, since it is clearly invertible. Indeed, to recover $\pi$ from $\hp\in\T_n$, write $\hp$ in cycle notation, replace the $\star$ with the entry in
$\{1,\dots,n\}$ that is missing, and turn the parentheses into brackets, thus recovering the one-line notation of $\pi$.

For $\hp\in\T_n$, let $\des(\hp)$ denote the number of descents of the sequence that we get by deleting the~$\star$ from the one-line notation of $\hp$. That is,
if $\hp=[a_1,\dots,a_j,\star,a_{j+1},\dots,a_{n-1}]$, then $\des(\hp)=|\{i:1\le i\le n-2,\, a_i>a_{i+1}\}|$.
We can now state a simpler formula for $N(\pi)$.

\begin{proposition}\label{prop:Npi2} Let $\pi\in\S_n$, $\hp=\theta(\pi)$. Then $N(\pi)$ is given by
$$N(\pi)=1+\des(\hp)+\epsilon(\hp),$$
where $$\epsilon(\hp)=\begin{cases}
1 & \mbox{if } \hp=[\star,1,\dots] \mbox{ or } \hp=[\dots,n,\star], \\
0 & \mbox{otherwise.}
\end{cases}$$
\end{proposition}

\begin{proof}
We compare this formula with the one in Theorem~\ref{th:Npi}. If $\pi(n)\notin\{1,n\}$, then $\epsilon(\hp)=0$, and $\des(\hp)=|A(\pi)|+1$ or $|A(\pi)|$ depending on whether $\pi(i+1)>\pi(j+1)$ or not, where
$i=\pi^{-1}(\pi(n)-1)$ and $j=\pi^{-1}(\pi(n)+1)$. If $\pi(n)\in\{1,n\}$, then $\des(\hp)=|A(\pi)|$ and $\epsilon(\hp)=\Delta(\pi)$.
In all cases, the two formulas for $N(\pi)$ are equivalent.
\end{proof}

For example, if $\pi=[8,9,2,3,6,4,1,5,7]$, then $\hp=[5,3,6,1,7,4,\star,9,2]$ has 4 descents, so $N(\pi)=1+4+0=5$.
If $\pi=[8,9,3,1,4,6,2,7,5]$, then $\hp=[4,7,1,6,\star,2,5,9,3]$ has 3 descents, so $N(\pi)=1+3+0=4$.
If $\pi=[3,4,2,1]$, then $\hp=[\star,1,4,2]$ has 1 descent, so $N(\pi)=1+1+1=3$.

If $\pi\in\S_n$, we have by definition that $N(\pi)=\min\{N:\pi\notin\F_n(\Sigma_N)\}=\min\{N:\pi\in\Al_n(\Sigma_N)\}$.
Proposition~\ref{th:aek_minshift} can now be recovered as a consequence of Proposition~\ref{prop:Npi2}.

\begin{corollary}
Let $n\ge3$. For every $\pi\in\S_n$, $N(\pi)\le n-1$.
\end{corollary}

\begin{proof}
Since $\hp$ has $n-1$ entries other than the $\star$, it can have at most $n-2$ descents. Besides, in the cases when $\epsilon(\hp)=1$,
we have necessarily $\des(\hp)\le n-3$. Thus, by Proposition~\ref{prop:Npi2}, $N(\pi)=1+\des(\hp)+\epsilon(\hp)\le n-1$ in any case.
\end{proof}

\ms

We define $\S_{n,N}=\{\pi\in\S_n:N(\pi)=N\}$. We are interested in the numbers $a_{n,N}=|\S_{n,N}|$. To avoid the trivial cases, we will assume that $n,N\ge2$.
From the definitions,
%$$\S_{n,N}=\F_n(\Sigma_{N-1})\setminus\F_n(\Sigma_{N})=\Al_n(\Sigma_{N})\setminus\Al_n(\Sigma_{N-1}),$$
$$\Al_n(\Sigma_M)=\bigcup_{N=2}^{M} \S_{n,N},\qquad \F_n(\Sigma_M)=\bigcup_{N=M+1}^{n-1} \S_{n,N}.$$
Since the sets $\S_{n,N}$ are disjoint, we have that $$|\Al_n(\Sigma_M)|=\sum_{N=2}^{M} a_{n,N},\qquad |\F_n(\Sigma_M)|=\sum_{N=M+1}^{n-1} a_{n,N}.$$
The first few values of $a_{n,N}$ are given in Table~\ref{tab:anN}. It follows from Lemma~\ref{lem:sym} that all the $a_{n,N}$ are even.

\begin{table}[hbt] \centering
\renewcommand{\baselinestretch}{5}
\begin{tabular}{r|r|r|r|r|r|r}
  $n \backslash N$ & 2 & 3 & 4 & 5 & 6 & 7 \\ \hline
  2 & 2 &  &  &  &  &  \\ \hline
  3 & 6 &  &  &  &  &  \\ \hline
  4 & 18 & 6 &  &  &  &  \\ \hline
  5 & 48 & 66 & 6 &  &  &  \\ \hline
  6 & 126 & 402 & 186 & 6 &  & \\ \hline
  7 & 306 & 2028 & 2232 & 468 & 6 &  \\ \hline
  8 & 738 & 8790 & 19426 & 10212 & 1098 & 6 \\
\end{tabular}

\caption{The numbers $a_{n,N}=|\{\pi\in\S_n:N(\pi)=N\}|$ for $n\le8$.}\label{tab:anN}
\end{table}

The next result shows that, independently of $n$, there are exactly six permutations of length $n$ that require the maximum number of symbols (i.e., $n-1$) in order to be realized. This settles a conjecture from~\cite{AEK}.
Given a permutation $\pi\in\S_n$, we will use $\pi^{rc}$ to denote the permutation such that $\pi^{rc}(i)=n+1-\pi(n+1-i)$ for $1\le i\le n$. If $\sigma$ is a marked cycle, then $\sigma^{rc}$ is defined similarly,
where if $\sigma(i)$ is the marked entry of $\pi$, then $\sigma^{rc}(n+1-i)$ is the marked entry of $\sigma^{rc}$.
It will be convenient to visualize $\pi\in\S_n$ as an $n\times n$ array with dots in positions $(i,\pi(i))$, for $1\le i\le n$.
The first coordinate refers to the row number, which increases from left to right,
and the second coordinate is the column number, which increases from bottom to top. Then, the array of $\pi^{rc}$ is obtained from the array of $\pi$ by
a 180-degree rotation. Of course, the array of $\pi^{-1}$ is obtained from the one of $\pi$ by reflecting it along the diagonal $y=x$. Notice also that
the cycle structure of $\pi$ is preserved in $\pi^{-1}$ and in $\pi^{rc}$. A marked cycle can be visualized in the same way, replacing the dot corresponding to the distinguished element with a $\star$.

[Should define $\sigma^{-1}$ when $\sigma$ is a marked cycle. Note that one can't just reverse the cycle notation of $\sigma$ with the star.]

\begin{proposition}\label{prop:6}
For every $n\ge3$, $a_{n,n-1}=6$.
\end{proposition}

\begin{proof}
First we show (this part was already proved in~\cite{AEK}) that $a_{n,n-1}\ge 6$, by giving six permutations in this set.
Let $m=\lceil n/2 \rceil$, and let $$\sigma=[n,n-1,\dots,m+1,\star,m,m-1,\dots,2], \quad \tau=[\star,1,n,n-1\dots,m+2,m,m-1,\dots,2]\in\T_n$$
(see Figure~\ref{fig:sigmatau}). Using Proposition~\ref{prop:Npi2}, it is easy to check that if $\hp\in\{\sigma,\sigma^{rc},\sigma^{-1},(\sigma^{-1})^{rc},\tau,\tau^{rc}\}$,
then $N(\pi)=n-1$, and that the six permutations in the set are different.
\begin{figure}[hbt] \centering
\epsfig{file=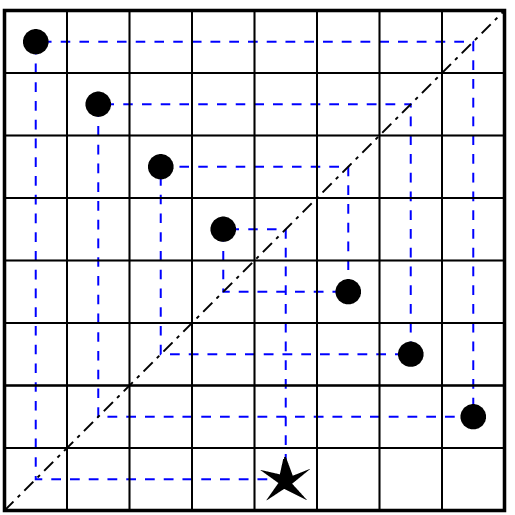,height=4cm} \hspace{1cm} \epsfig{file=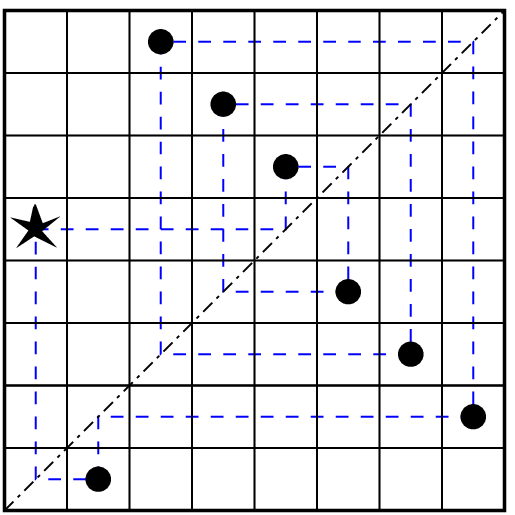,height=4cm} \caption{The arrays of $\sigma$ and $\tau$ for $n=8$, with dotted lines indicating the cycle structure.}\label{fig:sigmatau}
\end{figure}

Let us now show that there are no other permutations with $N(\pi)=n-1$.
We know by Proposition~\ref{prop:Npi2} that $N(\pi)=n-1$ can only happen if $\des(\hp)=n-2$, or if $\des(\hp)=n-3$ and $\epsilon(\hp)=1$.

{\it Case 1: $\des(\hp)=n-2$.} In this case, all the entries in $\hp$ other that the $\star$ must be in decreasing order. If the distinguished
entry is neither $\hp(1)$ nor $\hp(n)$, then the $\star$ must be replacing either $1$ or $n$; otherwise we would have that $\hp(1)=n$ and
$\hp(n)=1$, so $\hp$ would not be an $n$-cycle. It follows that in the array of $\hp$, the entry corresponding to the $\star$ is either in the
top or bottom row, or in the leftmost or rightmost column.

If the $\star$ is replacing $1$ (i.e, it is is the bottom row of the array), we claim that the only possible $n$-cycle in which the other entries are
in decreasing order is $\hp=\sigma$. Indeed, if we consider the cycle structure of $\hp=(1,\hp(1),\hp^2(1),\dots,\hp^{n-1}(1))$,
we see that $\hp(1)=n$ and $\hp^2(1)=\hp(n)=2$. Now, $\hp^i(1)\neq1$ for $3\le i\le n-1$, so the decreasing condition on the remaining entries
forces $\hp^3(1)=\hp(2)=n-1$, $\hp^4(1)=\hp(n-1)=3$, and so on.
A similar argument, considering that rotating the array 180 degrees preserves the cycle structure, shows that if the $\star$ is replacing $n$
(i.e, it is in the top row of the array), then necessarily $\hp=\sigma^{rc}$.

If the distinguished entry is $\hp(1)$ (i.e, it is in the leftmost column of the array), then a symmetric argument, reflecting the array along $y=x$),
shows that $\hp=\sigma^{-1}$. Similarly, if the distinguished entry is $\hp(n)$ (i.e, it is is the rightmost column of the array), then necessarily $\hp=(\sigma^{-1})^{rc}$.

{\it Case 2: $\des(\hp)=n-3$ and $\epsilon(\hp)=1$.} The second condition forces $\hp=[\star,1,\dots]$ or $\hp=[\dots,n,\star]$. Let us restrict to the
first case (the second one can be argued in a similar way if we rotate the array 180 degrees). We must have $\hp(3)>\hp(4)>\dots>\hp(n)$. We claim
that the only such $\hp$ that is also an $n$-cycle is $\hp=\tau$. Indeed, looking at the cycle structure
$\hp=(\hp^{-(n-1)}(1),\dots,\hp^{-1}(1),1)$, we see that $\hp^{-1}(1)=2$.
Now, $\hp^{-i}(1)\neq1$ for $2\le i\le n-1$, so the decreasing condition on the remaining entries forces $\hp^{-2}(1)=\hp^{-1}(2)=n$, $\hp^{-3}(1)=\hp^{-1}(n)=3$, $\hp^{-4}(1)=\hp^{-1}(3)=n-1$, and so on.
\end{proof}

\section{The number of allowed patterns of the binary shift}\label{sec:binary}

In the rest of the paper, we will assume for simplicity that $w_A(\pi)$ and $w_B(\pi)$ are defined taking $m=n-1$, so they are of the form $up^{n-1}x^\infty$, with $x=0$ or $x=N-1$ respectively.
The following variation of Lemma~\ref{lem:nk} will be useful later.

\begin{lemma}\label{lem:consec}
Let $w=up^{n-1}0^\infty\in\W_N$, where $|u|=k-1$ and $|p|=n-k$ for some $1\le k\le n-1$, and $p$ is primitive. If $\pi=\Pat(w,\Sigma,n)$ is defined, then $\pi(n)=\pi(k)-1$.
\end{lemma}

\begin{proof}
We have that $w_{[k,\infty)}=p^{n-1}0^\infty$ and $w_{[n,\infty)}=p^{n-2}0^\infty$. Since $p$ must have some nonzero entry (otherwise $\Pat(w,\Sigma,n)$ would not be defined),
it is clear that $w_{[n,\infty)}<w_{[k,\infty)}$. The same argument as in the proof of Lemma~\ref{lem:nk} can be used to show that there is no $1\le s\le n$ such that
$p^{n-2}0^\infty<w_{[s,\infty)}<p^{n-1}0^\infty$. It follows that $\pi(n)=\pi(k)-1$.
\end{proof}

For $n\ge2$, the set of patterns of length $n$ that are realized by the shift on two symbols is
$\Al_n(\Sigma_2)=\S_{n,2}$. In this section we enumerate these permutations. Recall that $a_{n,2}=|\S_{n,2}|$ and that $\psi_2(t)$ is the number of primitive binary words of length $t$.

\begin{theorem}\label{th:binary} For $n\ge2$,
$$a_{n,2}=\sum_{t=1}^{n-1}\psi_2(t)2^{n-t-1}.$$
\end{theorem}

\begin{proof}
Fix $n\ge 2$. We will construct a set $W\subset\W_2$ with the following four properties:
\begin{enumerate}\renewcommand{\labelenumi}{(\roman{enumi})}
\item for all $w\in W$, $\Pat(w,\Sigma_2,n)$ is defined,
\item for all $w,w'\in W$ with $w\neq w'$, we have that $\Pat(w,\Sigma_2,n)\neq\Pat(w',\Sigma_2,n)$,
\item for all $\pi\in\Al_n(\Sigma_2)$, there is a word $w\in W$ such that $\Pat(w,\Sigma_2,n)=\pi$,
\item $|W|=\sum_{t=1}^{n-1}\psi_2(t)2^{n-t-1}$.
\end{enumerate}
Properties (i)-(iii) imply that the map \bea W & \longrightarrow & \S_{n,2} \nn \\ w & \mapsto & \Pat(w,\Sigma_2,n) \label{eq:mapWSn2}\eea is a bijection. Thus, $a_{n,2}=|W|$ and the
result will follow from property (iv).

Let
$$W=\bigcup_{t=1}^{n-1}\{up^{n-1}x^\infty : u\in\{0,1\}^{n-t-1}, p\in\{0,1\}^t \mbox{ is a primitive word, and } x=\bar{p_t} \},$$
where we use the notation $\bar{0}=1$, $\bar{1}=0$. Given binary words $u,p$ of lengths $n-t-1$ and $t$ respectively,
where $p$ is primitive, and $x=\bar{p_t}$, we will denote $v(u,p)=up^{n-1}x^\infty$.

To see that $W$ satisfies (i), we have to show that for any $w\in W$ and any $1\le i<j\le n$, we have $w_{[i,\infty)}\neq w_{[j,\infty)}$. This is clear
because if $x=0$ (resp. $x=1$) both $w_{[i,\infty)}$ and $w_{[j,\infty)}$ end with $10^\infty$ (resp. $01^\infty$), with the last $1$ (resp. $0$) being
in different positions in $w_{[i,\infty)}$ and $w_{[j,\infty)}$.

Now we prove that $W$ satisfies (ii). Let $u$, $u'$ be binary words of lengths $n-t-1$, $n-t'-1$, respectively, and let $p$, $p'$ be primitive binary words
of lengths $t$, $t'$, respectively. Let $w=v(u,p)$ and $w'=v(u',p')$, and let $\pi=\Pat(w,\Sigma_2,n)$, $\pi'=\Pat(w',\Sigma_2,n)$.
We assume that $w\ne w'$, and want to show that $\pi\neq\pi'$. From $w\ne w'$ it follows that $u\ne u'$ or $p\ne p'$.

Corollary~\ref{cor:determined} for $N=2$ implies that if
$w_1w_2\dots w_{n-1}\neq w'_1w'_2\dots w'_{n-1}$, then $\Pat(w,\Sigma_2,n)\neq \Pat(w',\Sigma_2,n)$. In particular, if $t=t'$, then $up\neq u'p'$, so $\pi\neq\pi'$.

We are left with the case that $t\neq t'$%(we can assume without loss of generality that $t<t'$)
and $up=u'p'=w_1w_2\dots w_{n-1}$.
Let us first assume that $w_{n-1}=1$ (and so $p_t=p'_{t'}=1$). By Lemma~\ref{lem:consec} with $k=n-t$, we have that $\pi(n)=\pi(n-t)-1$, and similarly
$\pi'(n)=\pi'(n-t')-1$. If we had that $\pi=\pi'$, then $\pi(n)=\pi'(n)$ and so $\pi(n-t)=\pi'(n-t')=\pi(n-t')$. But $t\neq t'$, so this is a contradiction.
In the case $w_{n-1}=0$, an analogous argument to the proof of Lemma~\ref{lem:consec} implies that $w_{[n-t,\infty)}=p^{n-1}1^\infty<p^{n-2}1^\infty=w_{[n,\infty)}$ and there is no $s$
such that $w_{[s,\infty)}$ is strictly in between the two. Thus, $\pi(n)=\pi(n-t)+1$, and similarly $\pi'(n)=\pi'(n-t')+1$, so again $\pi\neq\pi'$.

To see that $W$ satisfies (iii) we use the construction from the proof of the upper bound in Theorem~\ref{th:Npi}.
Let $\pi\in\Al_n(\Sigma_2)$.
If $\pi(n-1)>\pi(n)$, let $w=w_A(\pi)=up^{n-1}0^\infty$. By Lemma~\ref{lem:nonzero}, $w_{n-1}=1$, so $w\in W$.
Similarly, if $\pi(n-1)<\pi(n)$, let $w=w_B(\pi)=up^{n-1}1^\infty$. By Lemma~\ref{lem:nonmax}, $w_{n-1}=0$, so $w\in W$. In both cases, $\Pat(w,\Sigma_2,n)=\pi$, so this construction is the inverse of the map~(\ref{eq:mapWSn2}).

To prove (iv), observe that the union in the definition of $W$ is a disjoint union. This is because the value of $t$ determines the position of the last entry
in $w$ that is not equal to $x$. For fixed $t$, there are $2^{n-t-1}$ choices for $u$ and $\psi_2(t)$ choices for $t$, so the formula follows.
\end{proof}

\begin{example} For $n=3$, we have $W=\{\ul{0}\,0\,0\,1^\infty\, \ul{0}\,1\,1\,0^\infty, \ul{1}\,0\,0\,1^\infty, \ul{1}\,1\,1\,0^\infty,\ 10\,10\,1^\infty, 01\,01\,0^\infty\}$,
where each word is written as $w=\ul{u}\,p\,p\,x^\infty$. The permutations corresponding to these words are
$\Al_3(\Sigma_2)=\{123,132,312,321,\ 231,213\}$, in the same order.

For $n=4$, we have
\bea \nn W=& \{\ul{00}\,0\,0\,0\,1^\infty, \ul{00}\,1\,1\,1\,0^\infty, \ul{01}\,0\,0\,0\,1^\infty, \ul{01}\,1\,1\,1\,0^\infty, \ul{10}\,0\,0\,0\,1^\infty, \ul{10}\,1\,1\,1\,0^\infty, \ul{11}\,0\,0\,0\,1^\infty, \ul{11}\,1\,1\,1\,0^\infty, \\
\nn & \ul{0}\,01\,01\,01\,0^\infty, \ul{0}\,10\,10\,10\,1^\infty, \ul{1}\,01\,01\,01\,0^\infty, \ul{1}\,10\,10\,10\,1^\infty, \\
\nn & 001\,001\,001\,0^\infty, 010\,010\,010\,1^\infty, 011\,011\,011\,0^\infty, 100\,100\,100\,1^\infty, 101\,101\,101\,0^\infty, 110\,110\,110\,1^\infty\},\eea
where each word is written as $w=\ul{u}\,p\,p\,p\,x^\infty$. The permutations corresponding to these words are
\bea \nn \Al_4(\Sigma_2)= & \{ 1234,1243,3412,1432,4123,2143,4312,4321, \\
\nn & 1342,1324,4231,4213, \\
\nn & 2341,2413,2431,3124,3142,3214 \}.\eea
\end{example}

\section{The number of allowed patterns of the shift on $N$ symbols}\label{sec:anyN}

In this section we give a formula for the numbers $a_{n,N}$, which count permutations that can be realized by the shift on $N$ symbols but not by the shift on $N-1$ symbols. This generalizes Theorem~\ref{th:binary}.
For any $N\ge 2$, define
\beq\label{eq:defomega}\Omega_{n,N}=\bigcup_{t=1}^{n-1}\{up^{n-1}0^\infty : u\in\{0,1,\dots,N-1\}^{n-t-1},\mbox{ and } p\in\{0,1,\dots,N-1\}^t \mbox{ is a primitive word}\}.\eeq

\begin{lemma}\label{lem:defined}
Let $w=up^{n-1}0^\infty\in\Omega_{n,N}$. Then $\Pat(w,\Sigma,n)$ is defined if and only if $p\neq 0$.
\end{lemma}

\begin{proof}
If $p=0$, then $w_{[n-1,\infty)}=w_{[n,\infty)}=0^\infty$, so $\Pat(w,\Sigma,n)$ is not defined.
Otherwise, since $p$ is primitive, it must have some nonzero entry. Let $p_\ell$ be the rightmost nonzero entry of $p$. Then, for any $1\le i<j\le n$, we have that $w_{[i,\infty)}\neq w_{[j,\infty)}$, because
both $w_{[i,\infty)}$ and $w_{[j,\infty)}$ end with $p_\ell0^\infty$, with the last nonzero entry being
in different positions in $w_{[i,\infty)}$ and $w_{[j,\infty)}$.
\end{proof}

Let
$\H'_{n,N}=\{\pi\in\S_{n,N}:\pi(n)=n\}$, and let $h_{n,N}=|\H'_{n,N}|$.
Recall the definition of $w_A(\pi)$ from Section~\ref{sec:upperbound}.
We have a map
$$\bt{ccc} $\S_{n,N}\setminus \H'_{n,N}$ & $\longrightarrow$ & $\Omega_{n,N}$ \\
 $\pi$ & $\mapsto$ & $w_A(\pi)$. \et $$
This is a one-to-one map, since $\pi$ can be recovered from $w=w_A(\pi)$ using that $\Pat(w,\Sigma,n)=\pi$.
Denote its range by $$\Gamma_{n,N}=\{w_A(\pi):\pi\in\S_{n,N}\setminus \H'_{n,N}\}.$$
%This is a disjoint union, because for any $\pi\in\S_n$, $w_A(\pi)$ ends with $0^\infty$ but $w_B(\pi)$ ends with $(N-1)^\infty$.
Defining $g_{n,N}=|\Gamma_{n,N}|$, it is clear that
$$g_{n,N}=|\S_{n,N}\setminus \H'_{n,N}|=a_{n,N}-h_{n,N}.$$
Our first goal is to find a formula for $g_{n,N}$.

\ms

Let $w=up^{n-1}0^\infty\in\Omega_{n,N}$. By Lemma~\ref{lem:defined}, $\Pat(w,\Sigma,n)$ is defined unless $p=0$. Assuming that $p\neq0$, let $\pi=\Pat(w,\Sigma,n)$. Clearly $\pi\in\Al_n(\Sigma_N)$, so $2\le N(\pi)\le N$.
Let $M=N(\pi)$, and let $j=N-M$. We can write $\Omega_{n,N}$ as a disjoint union:
\bea\label{eq:decompomega}\Omega_{n,N}&=&\bigcup_{j=0}^{N-2} \{w\in\Omega_{n,N} : \pi=\Pat(w,\Sigma,n) \mbox{ is defined and } N(\pi)=N-j\} \\ \nn &&\cup\quad\{u0^\infty : u\in\{0,1,\dots,N-1\}^{n-2}\} ,\eea
where the last term is the set of words $w$ for which $\Pat(w,\Sigma,n)$ is not defined.
Let us restrict to the case where $\pi=\Pat(w,\Sigma,n)$ is defined.

\begin{lemma}\label{lem:notn}
Let $w\in\Omega_{n,N}$. If $\pi=\Pat(w,\Sigma,n)$ is defined, then $\pi(n)\neq n$.
\end{lemma}

\begin{proof}
We have that $w_{[n,\infty)}=p^{n-2}0^\infty<p^{n-1}0^\infty$, so $\pi(n)$ is not the largest value in $\pi$.
\end{proof}

For $j=0$, we have the following result.

\begin{lemma}\label{lem:wwA}
Let $w\in\Omega_{n,N}$. If $\pi=\Pat(w,\Sigma,n)$ is defined and $N(\pi)=N$, then $w=w_A(\pi)$.
\end{lemma}

\begin{proof}
Write $w=up^{n-1}0^\infty$, and let $w'=w_A(\pi)=u'p'^{n-1}0^\infty\in\W_N$.
Since $\Pat(w',\Sigma,n)=\Pat(w,\Sigma,n)=\pi$, Corollary~\ref{cor:determined} implies that $u'p'=w'_1w'_2\dots w'_{n-1}=w_1w_2\dots w_{n-1}=up$, because these entries
are uniquely determined by $\pi$. Let $t=|p|$ and $t'=|p'|$. From the definition of $w_A$, we have that $\pi(n)=\pi(n-t')-1$, and by Lemma~\ref{lem:consec} with $k=n-t$, we have that $\pi(n)=\pi(n-t)-1$.
It follows that $\pi(n-t')=\pi(n-t)$, so $t=t'$. Therefore, $w=w'$.
\end{proof}

From Lemmas~\ref{lem:notn} and~\ref{lem:wwA} we have that
$$\{w\in\Omega_{n,N} : \pi=\Pat(w,\Sigma,n) \mbox{ is defined and } N(\pi)=N\}=\Gamma_{n,N}.$$
For general $j$, the situation is slightly more involved.

\begin{lemma}\label{lem:fewersymb}
Let $w\in\Omega_{n,N}$. If $\pi=\Pat(w,\Sigma,n)$ is defined and $N(\pi)=M=N-j$, then the entries $w_1w_2\dots w_{n-1}$ satisfy
$$(w_{i_1},w_{i_2},\dots,w_{i_{n-1}})=(v_{i_1}+c_1,v_{i_2}+c_2,\dots,v_{i_{n-1}}+c_{n-1})$$
for some $0\le c_1\le c_2\le \dots \le c_{n-1}\le j$, where $i_1,\dots,i_{n-1}$ are the subindices in equation~(\ref{eq:increasing_n}) from left to right, and $v=w_A(\pi)\in\Gamma_{n,M}$.
\end{lemma}

\begin{proof}
First note that by Lemma~\ref{lem:notn}, $\pi(n)\neq n$, so $v\in\Gamma_{n,M}$.
We know by Theorem~\ref{th:Npi} that $M=1+|A(\pi)|+\Delta(\pi)$.
By definition of $i_1,i_2,\dots,i_{n-1}$, equation~(\ref{eq:increasing_n}) can be written as \beq\label{eq:wi} w_{i_1}\le w_{i_2}\le \dots \le w_{i_{n-1}}. \eeq
By Lemma~\ref{lem:required}, the entries of $w$ have to fulfill~(\ref{eq:wi}), plus $M-1$ required strict inequalities. Similarly, since $\Pat(v,\Sigma_M,n)=\pi$, the entries of $v$ satisfy
$v_{i_1}\le v_{i_2}\le \dots \le v_{i_{n-1}}$ and the same $M-1$ required strict inequalities. In fact, by Corollary~\ref{cor:determined}, the entries $v_1v_2\dots v_{n-1}$ are uniquely determined by
these inequalities and the fact that $v\in\W_M$.

It is now easy to see that $w\in\Omega_{n,N}$ satisfies~(\ref{eq:wi}) plus the $M-1$ required strict inequalities if and only if the values $w_1w_2\dots w_{n-1}$ satisfy
$$(w_{i_1},w_{i_2},\dots,w_{i_{n-1}})=(v_{i_1}+c_1,v_{i_2}+c_2,\dots,v_{i_{n-1}}+c_{n-1})$$
for some $0\le c_1\le c_2\le \dots \le c_{n-1}\le j$.
\end{proof}

Note that for $j=0$, Lemma~\ref{lem:fewersymb} says that the values $w_1w_2\dots w_{n-1}$ are determined by $\pi$, which is just Corollary~\ref{cor:determined}.
For general $j$, there are $\binom{n+j-1}{j}$ choices for the vector $(c_1,c_2,\dots,c_{n-1})$.
Let $C_{n,j}=\{(c_1,c_2,\dots,c_{n-1}):0\le c_1\le c_2\le \dots \le c_{n-1}\le j\}$.

\begin{lemma}\label{lem:bij} For $0\le j\le N-2$, the map $$\begin{array}{ccc} \{w\in\Omega_{n,N} : \pi=\Pat(w,\Sigma,n) \mbox{ is defined and } N(\pi)=N-j\}&\longrightarrow&\Gamma_{n,N-j}\times C_{n,j} \\
\nn  w & \mapsto &(v,(c_1,c_2,\dots,c_{n-1}))\end{array}$$ defined by Lemma~\ref{lem:fewersymb} is a bijection.
\end{lemma}

\begin{proof}
The map is clearly well defined, since $v=w_A(\pi)\in\Gamma_{n,N-j}$ and $(c_1,\dots,c_{n-1})\in C_{n,j}$ are uniquely determined by $w$.
To see that it is invertible, take $(v,(c_1,c_2,\dots,c_{n-1}))\in\Gamma_{n,N-j}\times C_{n,j}$. Then, $\pi=\Pat(v,\Sigma,n)$ determines the sequence $i_1,i_2,\dots,i_{n-1}$ of subindices in equation~(\ref{eq:increasing_n}),
so the first $n-1$ entries in $w$ must be $(w_{i_1},w_{i_2},\dots,w_{i_{n-1}})=(v_{i_1}+c_1,v_{i_2}+c_2,\dots,v_{i_{n-1}}+c_{n-1})$. By Lemma~\ref{lem:consec}, the unique word in $w\in\Omega_{n,N}$ with these entries such that
$\Pat(w,\Sigma,n)=\pi$ is $w=w_1w_2\dots w_{k-1}(w_kw_{k+1}\dots w_{n-1})^{n-1}0^\infty$, where $k=\pi^{-1}(\pi(n)+1)$.
\end{proof}

\begin{proposition}\label{prop:g} The numbers $g_{n,N}$ satisfy
$$\sum_{t=1}^{n-1}\psi_N(t)N^{n-t-1}=\sum_{j=0}^{N-2}\binom{n+j-1}{j}g_{n,N-j}+N^{n-2}.$$
\end{proposition}

\begin{proof}
The formula is obtained by taking set cardinalities in equation~(\ref{eq:decompomega}).
For the left hand side we use that $|\Omega_{n,N}|=\sum_{t=1}^{n-1}\psi_N(t)N^{n-t-1}$, which follows from the definition~(\ref{eq:defomega}) of $\Omega_{n,N}$, noticing that it is a disjoint union.
By Lemma~\ref{lem:bij}, $$|\{w\in\Omega_{n,N} : \pi=\Pat(w,\Sigma,n) \mbox{ is defined and } N(\pi)=N-j\}|=|\Gamma_{n,N-j}||C_{n,j}|=g_{n,N-j}\binom{n+j-1}{j}$$ for $0\le j\le N-2$. Finally,
it is clear that $|\{u0^\infty : u\in\{0,1,\dots,N-1\}^{n-2}\}|=N^{n-2}$.
\end{proof}

From Proposition~\ref{prop:g} we obtain the following recurrence for $g_{n,N}$, for $N\ge 2$.
\beq\label{eq:recg} g_{n,N}=\sum_{t=1}^{n-1}\psi_N(t)N^{n-t-1}-N^{n-2}-\sum_{j=1}^{N-2}\binom{n+j-1}{j}g_{n,N-j}.\eeq
%Its initial term is $g_{n,2}=\sum_{t=1}^{n-1}\psi_2(t)2^{n-t-1}$, which agrees with Theorem~\ref{th:binary}.
To obtain an expression for the general term $g_{n,N}$, we use the following lemma.

\begin{lemma}\label{lem:solverec}
Let $n\ge1$ and $\{b_N\}_{N\ge 2}$ be fixed. If the sequence $\{r_N\}_{N\ge 2}$ satisfies the recurrence
$$r_N=b_N-\sum_{j=1}^{N-2}\binom{n+j-1}{j}r_{N-j},$$
then
$$r_N=\sum_{i=0}^{N-2}(-1)^i\binom{n}{i}b_{N-i}$$
for all $N\ge 2$.
\end{lemma}

\begin{proof}
We proceed by induction on $N$. Assuming that $r_M=\sum_{i=0}^{M-2}(-1)^i\binom{n}{i}b_{M-i}$ for $2\le M\le N-1$, we have that
\bea\nn r_N&=&b_N-\sum_{j=1}^{N-2}\binom{n+j-1}{j}\left(\sum_{i=0}^{N-j-2}(-1)^i\binom{n}{i}b_{N-j-i}\right)\\
\nn &=&b_N-\sum_{k=1}^{N-2}\sum_{j=1}^{k}\binom{n+j-1}{j}(-1)^{k-j}\binom{n}{k-j}b_{N-k}\\
\nn &=&b_N-\sum_{k=1}^{N-2}(-1)^{k-1}\binom{n}{k}b_{N-k}\ =\ \sum_{k=0}^{N-2}(-1)^k\binom{n}{k}b_{N-k}, \eea
where in the second step we have made the substitution $k=i+j$,
and in the third step we have used the identity $\sum_{j=1}^k (-1)^{j-1}\binom{n+j-1}{j}\binom{n}{k-j}=\binom{n}{k}$.
\end{proof}

Substituting $b_N=\sum_{t=1}^{n-1}\psi_N(t)N^{n-t-1}-N^{n-2}$ in Lemma~\ref{lem:solverec}, we can solve recurrence~(\ref{eq:recg}) for $g_{n,N}$.

\begin{corollary}\label{cor:g}
For any $n,N\ge2$,
$$g_{n,N}=\sum_{i=0}^{N-2}(-1)^i\binom{n}{i}\left(\sum_{t=1}^{n-1}\psi_{N-i}(t)(N-i)^{n-t-1}-(N-i)^{n-2}\right).$$
\end{corollary}

\bs

Recall that $g_{n,N}=a_{n,N}-h_{n,N}$, and that our main goal is to compute $a_{n,N}$, the number of permutations of length $n$ that require $N$ symbols
in the alphabet in order to be realized by a shift. The next step will be to find a formula for $h_{n,N}=|\{\pi\in\S_{n,N}:\pi(n)=n\}|$, using similar ideas to the ones
we used for $g_{n,N}$.
Let $\H_{n,N}=\{\pi\in\S_{n,N}:\pi(n)=1\}$. By Lemma~\ref{lem:sym}, $h_{n,N}=|\H_{n,N}|$. We now define analogues of $\Omega_{n,N}$ and $\Gamma_{n,N}$ corresponding
to permutations ending with a 1. Let
$$\Theta_{n,N}=\{w\in\Omega_{n,N} : \pi=\Pat(w,\Sigma,n) \mbox{ is defined and has } \pi(n)=1\},$$
$$\Lambda_{n,N}=\{w_A(\pi):\pi\in\H_{n,N}\}.$$
Since $w_A$ is one-to-one, it is clear that $|\Lambda_{n,N}|=h_{n,N}$.

\begin{proposition}\label{prop:Theta}
For any $n,N\ge2$, $|\Theta_{n,N}|=(N-1)N^{n-2}$.
\end{proposition}

\begin{proof}
We will give a bijection between the set $$R_{n,N}=\{w_1w_2\dots w_{n-1} : w_i\in\{0,1,\dots,N-1\} \mbox{ for } 1\le i\le N-2, \mbox{ and } w_{n-1}\in\{1,2,\dots,N-1\}\}$$
and $\Theta_{n,N}$.

Given $w_1\dots w_{n-1}\in R_{n,N}$, we want to construct a word $w=up^{n-1}0^\infty\in\Theta_{n,N}$ where $u=w_1w_2\dots w_{k-1}$ and $p=w_k w_{k+1}\dots w_{n-1}$ for some $1\le k\le n-1$.
Let $k$ be the index such that $w_k w_{k+1}\dots w_{n-1}0^\infty$ is smallest (in lexicographic order), with $1\le k\le n-1$. Note that these $n-1$ words
are all different because the last nonzero entry of $w_k w_{k+1}\dots w_{n-1}0^\infty$ is $w_{n-1}$, which is in position $n-k$.
We claim that with this choice of $k$, the word $p=w_k w_{k+1}\dots w_{n-1}$ is primitive. Otherwise, if $p=q^m$ for some $m>1$, then $q^m0^\infty>q0^\infty$, contradicting the choice of $k$.
Since $p\neq0$, Lemma~\ref{lem:defined} guarantees that $\pi=\Pat(w,\Sigma,n)$ is defined. Moreover, $\pi(n)=1$ because $w_{[n,\infty)}=p^{n-2}0^\infty<w_{[k,\infty)}=p^{n-1}0^\infty$, and
by the choice of $k$, $w_{[k,\infty)}<w_{[i,\infty)}$ for all $1\le i\le n-1$ with $i\neq k$. Thus $w\in\Theta_{n,N}$.

This map is clearly one-to one, because $w_1\dots w_{n-1}$ can be recovered by taking the first $n-1$ entries in $w$. To show that it is onto, we have to show that if
$w\in\Theta_{n,N}$, then $w_{n-1}\neq0$. Assume for contradiction that $w_{n-1}=0$. Then, $w_{[n-1,\infty)}=0w_{[n,\infty)}$, so $w_{[n-1,\infty)}\le w_{[n,\infty)}$, which contradicts the fact that $\pi=\Pat(w,\Sigma,n)$ has $\pi(n)=1$.

Clearly $|R_{n,N}|=(N-1)N^{n-2}$, so the proposition is proved.
\end{proof}

Similarly to the decomposition of $\Omega_{n,N}$ in equation~(\ref{eq:decompomega}), we can decompose $\Theta_{n,N}$ as
\beq\label{eq:decomptheta}\Theta_{n,N}=\bigcup_{j=0}^{N-2} \{w\in\Theta_{n,N} : N(\Pat(w,\Sigma,n))=N-j\}.\eeq
For $j=0$, Lemma~\ref{lem:wwA} implies that
$$\{w\in\Theta_{n,N} : N(\Pat(w,\Sigma,n))=N\}=\Lambda_{n,N}.$$
In general, given $w\in\Theta_{n,N}$, let $\pi=\Pat(w,\Sigma,n)$, and assume that $N(\pi)=N-j$. Then, $\pi\in\H_{n,N-j}$, and the word $v=w_A(\pi)$ given by
Lemma~\ref{lem:fewersymb} is in $\Lambda_{n,N-j}$. Restricting to permutations with $\pi(n)=1$, we obtain the following analogue of Lemma~\ref{lem:bij}.

\begin{lemma}\label{lem:bij2} For $0\le j\le N-2$, the map $$\begin{array}{ccc} \{w\in\Theta_{n,N} : N(\Pat(w,\Sigma,n))=N-j\}&\longrightarrow&\Lambda_{n,N-j}\times C_{n,j} \\
\nn  w & \mapsto &(v,(c_1,c_2,\dots,c_{n-1}))\end{array}$$ defined by Lemma~\ref{lem:fewersymb} is a bijection.
\end{lemma}

\begin{proposition}\label{prop:h} The numbers $h_{n,N}$ satisfy
$$(N-1)N^{n-2}=\sum_{j=0}^{N-2}\binom{n+j-1}{j}h_{n,N-j}.$$
\end{proposition}

\begin{proof}
The result follows by taking set cardinalities in equation~(\ref{eq:decomptheta}).
The left hand side is given by Proposition~\ref{prop:Theta}.
For the right hand side, we use Lemma~\ref{lem:bij2} to deduce that $$|\{w\in\Theta_{n,N} : N(\Pat(w,\Sigma,n))=N-j\}|=|\Lambda_{n,N-j}||C_{n,j}|=h_{n,N-j}\binom{n+j-1}{j}$$ for $0\le j\le N-2$.
\end{proof}

From Proposition~\ref{prop:h} we obtain the recurrence
$$h_{n,N}=(N-1)N^{n-2}-\sum_{j=1}^{N-2}\binom{n+j-1}{j}h_{n,N-j},$$
which can be solved using Lemma~\ref{lem:solverec} to get the following expression for the general term $h_{n,N}$.

\begin{corollary}\label{cor:h}
For any $n,N\ge2$,
$$h_{n,N}=\sum_{i=0}^{N-2}(-1)^i\binom{n}{i}(N-i-1)(N-i)^{n-2}.$$
\end{corollary}

Combining Corollaries~\ref{cor:g} and~\ref{cor:h}, and using the fact that $a_{n,N}=g_{n,N}+h_{n,N}$, we obtain an expression for the numbers $a_{n,N}=|\{\pi\in\S_n:N(\pi)=N\}|$.

\begin{theorem}\label{th:anN}
For any $n,N\ge2$,
\beq\label{eq:anN} a_{n,N}=\sum_{i=0}^{N-2}(-1)^i\binom{n}{i}\left((N-i-2)(N-i)^{n-2}+\sum_{t=1}^{n-1}\psi_{N-i}(t)(N-i)^{n-t-1}\right).\eeq
\end{theorem}

An asymptotic analysis of this formula shows that, as $n$ tends to infinity, $a_{n,N}\sim nN^{n-1}$. For $N=2$, Theorem \ref{th:anN} gives $a_{n,2}=\sum_{t=1}^{n-1}\psi_2(t)2^{n-t-1}$, which agrees with Theorem~\ref{th:binary}.

\section{Conjectures and further work}
\label{sec:conj}

Some general open questions that deserve further study have already been discussed at the end of Section~\ref{sec:forbdef}. Here we mention two curious conjectures that came up while studying forbidden patterns
of shift systems. They are derived from experimental evidence, and it would be interesting to find combinatorial proofs.

\subsection{An appearance of the Eulerian numbers}

Using Proposition~\ref{prop:Npi2}, the problem of enumerating $\S_{n,N}$ can be formulated in terms of counting marked cycles $\hp\in\T_n$ with respect to $\des(\hp)$ and $\epsilon(\hp)$.
More precisely, \beq\label{eq:sumdeseps} a_{n,N}=|\{\hp\in\T_n\,:\,\des(\hp)=N-1 \mbox{ and }\epsilon(\hp)=0\}|+|\{\hp\in\T_n\,:\,\des(\hp)=N-2 \mbox{ and }\epsilon(\hp)=1\}|. \eeq
Let us take a closer look at rightmost summand. Marked cycles with $\epsilon(\hp)=1$ can be separated into two disjoint sets $\E_n$ and $\E'_n$,
depending on whether they are of the form $\hp=[\star,1,\dots]$ or $\hp=[\dots,n,\star]$, respectively. The map $\hp\mapsto\hp^{rc}$ (see the definition above Proposition~\ref{prop:6}) is a bijection between
$\E_n$ and $\E'_n$, which preserves the number of descents. Thus, for any $0\le k\le n-3$,
$$|\{\hp\in\T_n\,:\,\des(\hp)=k \mbox{ and }\epsilon(\hp)=1\}|=2\ |\{\hp\in\E_n\,:\,\des(\hp)=k\}|.$$

It will be convenient to define the set $\T^0_n$ of $n$-cycles where one entry has been replaced with $0$. The set $\T^0_n$ is
essentially the same as $\T_n$, with the only difference that the $\star$ symbol in each element is replaced with a $0$. The reason for this change is that
now we define a descent of $\sigma=[\sigma(1),\dots,\sigma(n)]\in\T^0_n$ to be a position $i$ such that $\sigma(i)>\sigma(i+1)$, so the
$0$ entry is no longer skipped in the computation of the descent set $D(\sigma)$ and the number of descents $\des(\sigma)$.

\begin{lemma}\label{lem:Phi}
For $n\ge3$, the map
$$\bt{cccc} $\Phi:$ & $\E_n$ & $\longrightarrow$ & $\T^0_{n-2}$ \\ & $[\star,1,\hp(3),\hp(4),\dots,\hp(n)]$ & $\mapsto$ & $[\hp(3)-2,\hp(4)-2,\dots,\hp(n)-2]$\et$$
is well-defined and it is a bijection.
\end{lemma}

\begin{proof}
First we show that if $\hp\in\E_n$, then $\Phi(\hp)\in\T^0_{n-2}$. Let $j\in\{1,\dots,n\}$ be the missing entry in $\hp$, that is, the entry that $\star$ is replacing.
Clearly $j\neq 2$, otherwise $\hp$ would not be an $n$-cycle. In particular, there is some $3\le i\le n$ such that $\hp(i)=2$, which produces the 0 entry in $\Phi(\hp)$.
Besides, if we replace this 0 entry in $\Phi(\hp)$ with $j-2$, we obtain an $(n-2)$-cycle. This is indeed the same $(n-2)$-cycle that we get starting from the cycle notation of $\hp$,
replacing the $\star$ with $j$, deleting the entries $1$ and $2$, and shifting the values of the remaining entries down by 2. This proves that $\Phi(\hp)\in\T^0_{n-2}$.
%To show that $\Phi(\hp)\in\T^0_{n-2}$, it suffices to show that if the $0$ in $\Phi(\hp)$ is replaced with the element in
%$\{1,2,\dots,n-2\}$ that is missing from the entries in $\Phi(\hp)$, then the remaining permutation is an $n$-cycle. The reason for this is that if we write $\hat\pi$ in cycle notation,
%$\hp=(2,1,\star,\dots)$, then
%The entries of $\Phi(\hp)$ are the elements $\{1,2,\dots,n-2\}$ with one entry replaced with a 0.

The map is clearly invertible, and the fact that $[\hat\sigma(1),\hat\sigma(2),\dots,\hat\sigma(n-2)]\in\T^0_{n-2}$ guarantees that $[\star,1,\hat\sigma(1)+2,\hat\sigma(2)+2,\dots,\hat\sigma(n-2)+2]$ is a marked cycle.
\end{proof}

For example, if $\hp=[\star,1,5,7,6,2,3]$, then $\Phi(\hp)=[3,5,4,0,1]$. Note that from the cycle notation of $\hp=(2,1,\star,7,3,5,6)$, after replacing the $\star$ with a $4$, deleting $2$ and $1$, and shifting the remaining entries down
by 2, we obtain $(2,5,1,3,4)$, which is the cycle notation of $\Phi(\hp)$ with the 0 replaced with a 2.

It is clear from the definition of $\Phi$ that it preserves the number of descents. As a consequence of this property,
$$|\{\hp\in\E_n\,:\,\des(\hp)=k\}|=|\{\hp\in\T^0_{n-2}\,:\,\des(\hp)=k\}|,$$
so we have reduced the computation of the rightmost summand of~(\ref{eq:sumdeseps}) to studying the distribution of the number of descents in $\T^0_{n-2}$.
Experimental evidence suggests that this is a very well-known distribution:

\begin{conj} For any $n$ and any subset $D\subseteq\{1,2,\dots,n-1\}$,
$$|\{\sigma\in\T^0_n\,:\,D(\sigma)=D\}|=|\{\pi\in\S_n\,:\,D(\pi)=D\}|.$$
In particular, the statistic $\des$ has the same distribution in $\T^0_n$ as in $\S_n$, i.e,
$$\sum_{\sigma\in\T^0_n} x^{\des(\sigma)+1}=A_n(x),$$
the $n$-th Eulerian polynomial.
\end{conj}

We have checked this conjecture by computer for $n$ up to $9$.

\subsection{Divisibility properties}

Even though Theorem~\ref{th:anN} gives an explicit formula for the numbers $a_{n,N}$, some properties of these numbers are not apparent from the formula.
For example, it is not trivial to derive Proposition~\ref{prop:6} from equation~(\ref{eq:anN}). A striking property of the entries of Table~\ref{tab:anN} (and of other values of $a_{n,N}$ computed for larger $n$
using Theorem~\ref{th:anN}) is that they are all divisible by $6$. We conjecture that this is always the case.

\begin{conj}
For every $n,N\ge3$, $a_{n,N}$ is divisible by $6$.
\end{conj}

\section*{Acknowledgements}

The author is grateful to Jos\'e Mar\'{\i}a Amig\'o for useful comments on the
paper.

\end{document}